\documentclass[12pt,a4paper]{article}

\usepackage{latexsym}
\usepackage{amssymb}

\newcommand{\N}{\mathbb N}
\newcommand{\R}{\mathbb R}
\newcommand{\C}{\mathbb C}
\newcommand{\Z}{\mathbb Z}

\newcommand{\1}{\mathbf 1}

\newcommand{\Hom}{\mathrm{Hom}}
\newcommand{\Aut}{\mathrm{Aut}}
\newcommand{\tr}{\mathrm{tr}}
\newcommand{\Tr}{\mathrm{Tr}}
\newcommand{\ad}{\mathrm{ad}}
\newcommand{\id}{\mathrm{id}}
\renewcommand{\char}{\mathrm{char}}
\newcommand{\Alt}{\mathrm{Alt}}
\newcommand{\Sym}{\mathrm{Sym}}
\newcommand{\im}{\mathrm{im}}
\newcommand{\GL}{\mathrm{GL}}
\newcommand{\Mat}{\mathrm{Mat}}
\newcommand{\Bia}{\mathrm{Bia}}
\newcommand{\diag}{\mathrm{diag}}
\newcommand{\rk}{\mathrm{rk}}
\newcommand{\sgn}{\mathrm{sgn}}
\newcommand{\Deg}{\mathrm{Deg}}
\newcommand{\CO}{\mathrm{CO}}
\renewcommand{\O}{\mathrm{O}}
\newcommand{\SO}{\mathrm{SO}}
\newcommand{\Met}{\mathrm{Met}}
\newcommand{\Aff}{\mathrm{Aff}}

\newcommand{\sub}{\subseteq}
\newcommand{\Lie}[2]{[#1,#2]}
\newcommand{\iprod}[2]{\langle#1,#2\rangle}
\renewcommand{\phi}{\varphi}

\title{Bianchi's classification of $3$-dimensional Lie algebras revisited}

\author{Manuel Glas\thanks{Robert Bosch GmbH, Schwieberdingen} 
\and Panagiotis Konstantis\thanks{Fachbereich Mathematik und Informatik der Universit\"at Marburg} 
\and Achim Krause\thanks{Mathematisches Institut der Universit\"at T\"ubingen} 
\and Frank Loose\thanks{Mathematisches Institut der Universit\"at T\"ubingen}} 

\newtheorem{remark}{Remark}[section]
\newtheorem{lemma}[remark]{Lemma}
\newtheorem{corollary}[remark]{Corollary}
\newtheorem{theorem}[remark]{Theorem}
\newtheorem{proposition}[remark]{Proposition}

\newenvironment{proof}{Proof.}{\hfill$\Box$}

\begin{document}

\maketitle

\begin{abstract}
We present Bianchi's proof on the classification of real (and complex) $3$-dimensional Lie algebras in a coordinate free version from a strictly representation theoretic point of view. Nearby we also compute the automorphism groups and from this the orbit dimensions of the corresponding orbits in the algebraic variety $X\sub\Lambda^2V^*\otimes V$ describing all Lie brackets on a fixed vector space $V$ of dimension $3$. Moreover we clarify which orbits lie in the closure of a given orbit and therefore the topology on the orbit space $X/G$ with $G=\Aut(V)$.
\end{abstract}

%\tableofcontents

%\newpage

\setcounter{section}{-1}

\section{Introduction}

\noindent
It is now more than $115$ years ago when L.\ Bianchi presented his list of $3$-dimensional real (and complex) Lie algebras together with a proof that every $3$-dimensional Lie algebra is isomorphic to one, and only one, Lie algebra of his list (cf.\ \cite{Bianchi_1898,Bianchi_2001}). Due to him these isomorphism classes are noted even today by $\Bia(I)$, $\Bia(II)$, $\Bia(IV)$, $\Bia(V)$, $\Bia(VI_h)$ ($h\leq 0)$, $\Bia(VII_0)$ ($h\geq 0$), $\Bia(VIII)$ and $\Bia(IX)$ and we will do so in this paper.

The Bianchi list plays a vital role also in nowadays study of $3$-dimensional differential topology and in general relativity. The reason is perhaps first of all, because W.\ Thurston defined a $3$-dimensional geometry in the spirit of F.\ Klein as a pair $(M,G)$, where $M$ is $3$-dimensional simply connected manifold and $G$ is a Lie group, acting effectively and transitively on $M$ with compact stabilizer groups (cf.\ \cite{Thurston}). Moreover Thurston requires the existence of a discrete subgroup $\Gamma\sub G$, acting freely and cocompactly on $M$ and moreover a maximality for $G$ w.r.t.\ these properties. He then finds exactly $8$ model geometries. Using these he formulated his famous geometrisation program in which every closed §3§-manifold is built up from these $8$ model geometries using well defined procedures in differential topology. As we know, G.\ Perelman was able to prove this geometrisation conjecture in 2003.

For applications in cosmology, however, it is much more natural to resign on the existence of compact quotients and, if at all, demand a minimality for $G$ with these properties. P.\ Konstantis only recently classified all geometries in this more general setting (see \cite{Konstantis}) and among these geometries are, of course, the $3$-dimensional simply connected Lie groups, acting on itself, i.e., $M=G$ in this case, by left translation. Of course, these Lie groups are in a 1:1-correspondence to the $3$-dimensional Lie algebras. It is quite natural that these geometries $(M,G)$ appear as models for universa $L=I\times M$ ($I$ being an open interval in $\R$), which are in space directions not necessarily isotropic, as in the standard models, but at least acting transitively on the space leaves and therefore representing a spacelike homogeneity. The Lorentz metric $h$ is given here by $h=-dt^2+g(t)$, where $t$ is the coordinate of $I$ and $g\colon I\to\Met(M)$ is a curve of $G$-invariant Riemannian metrics on $M$. A good knowledge of the Bianchi algebras, together with the curvature formulas for the left-invariant metrics on the corresponding Lie groups, is therefore highly important.

In this paper we present a proof which takes a consequently representation theoretic point of view. Nearby we are able to compute also the automorphism groups of all Lie algebras and therefore in particular the orbit dimension of the corresponding orbits in the Lie variety $X\sub\Lambda^2V^*\otimes V$ of all Lie brackets $C=\Lie{\cdot}{\cdot}$ on a fixed vector space $V$ of dimension $3$. Finally, the whole orbit structure on $X$, i.e., the topology of the orbit space $X/G$, becomes clear.

\section{Classification of Bianchi A classes }

%Abschnitt 1: Bianchi-Klassifikation der 3-dimensionalen Liealgebren, Klasse A%

{\bf 1.1.} Let $V$ be a fixed vector space (over some field $K$ which is usually the field of real numbers $\R$ or of complex numbers $\C$) of dimension $n\in\N_0$.
A Lie bracket $C=\Lie{\cdot}{\cdot}$ on $V$ is first of all a bilinear map $C\colon V\times V\to V$ being skew symmetric, i.e.\
\[
 C(a,b)=-C(b,a),\qquad\mbox{for all $a,b\in V$}.
\]
Therefore $C$ induces a unique linear map $\tilde C\colon\Lambda^2V\to V$ with $\tilde C\circ\pi=C$ (and is determined by it), where $\pi\colon V\times V\to\Lambda^2V$ is the universal skew symmetric bilinear map $(v_1,v_2)\mapsto v_1\wedge v_2$. Therefore we look at $C$ as an element in
\[
 \Hom(\Lambda^2V,V)\cong(\Lambda^2V)^*\otimes V\cong\Lambda^2V^*\otimes V,
\]
$C\in\Lambda^2V^*\otimes V$. Additionally a Lie bracket fulfills the Jacobian identity. We come back to this later on.
\bigskip

\noindent
{\bf 1.2.} When are two Lie structures $C_1$ and $C_2$ on a vector space $V$ are equivalent? That means that there exists a Lie isomorphism between $(V,C_1)$ and $(V,C_2)$, i.e.\ an automorphism $\varphi\colon V\to V$ satisfying
\[
 \varphi\circ C_1(a,b)=C_2(\varphi a,\varphi b),
\]
for all $a,b\in V$. It is enough to consider this question for a given number $n\in\N_0$ on a fixed vector space $V$ of dimension $n$ since all vector spaces of dimension $ n$ are isomorphic (as vector spaces), of course. 

We consider therefore now the (algebraic) group $G:=\Aut(V)$, which we consider as a real resp.\ complex Lie group in case of $K=\R$ resp.\ $K=\C$. It acts in a natural way on $V$ and therefore also on all tensor products $T^{r,s}V=V^{\otimes r}\otimes(V^*)^{\otimes s}$ as well as on all outer products $\Lambda^rV$, on the symmetric products $S^kV$ ($r,s,k\in\N_0$) or on tensor products of these, e.g.\ on the vector space 
\[
 W=\Lambda^3V^*\otimes S^2V.
\]
All these are {\sl representations} of $G$, i.e.: there is a homomorphism $\rho\colon G\to\Aut(W)$.

\begin{remark} Let $V$ be a $K$ vector space. Two Lie structures $C_1,C_2\in\Lambda^2V^*\otimes V$ are equivalent if and only if they are in the same $G$ orbit in $W:=\Lambda^2V^*\otimes V$.
\end{remark}

\noindent\begin{proof}
The induced action of $G$ on $W$ is given by 
\[
 (\varphi.C)(a,b)=\varphi\circ C(\varphi^{-1}a,\varphi^{-1}b),
\]
for $\varphi\in G$ and $a,b\in V$. Therefore $\varphi$ is an isomorphism from $(V,C_1)$ to $(V,C_2)$ if and only if $\varphi^{-1}.C_1=C_2$, i.e., $C_1$ and $C_2$ are in the same orbit of $G$ in $W$.
\end{proof}
\bigskip

\noindent
{\bf 1.3.} In order to classify all Lie algebras of a fixed dimension $n\in\N_0$ one has therefore to do the following two things:

\begin{description}
 \item[(i)] The Jacobi condition
\[
 C(C(a,b),c)+C(C(b,c),a)+C(C(c,a),b)=0,
\]
for all $a,b,c\in V$, yields an algebraic subvariety $X\sub\Lambda^2V^*\otimes V\cong K^N$ (with $N={n\choose 2}n$), which is $G$ invariant, i.e., for $C\in X$ and $\phi\in G$  we have $\phi.C\in G$ as well, since
\begin{eqnarray*}
 \lefteqn{\phi.C(\phi.C(a,b),c)+\phi.C(\phi.C(b,c),a)+\phi.C(\phi.C(c,a),b)} \\
  & = & \phi.C(\phi\circ C(\phi^{-1}a,\phi^{-1}b),c)+\cdots \\
  & = & \phi\circ C(\phi^{-1}\circ\phi\circ C(\phi^{-1}a,\phi^{-1}b),\phi^{-1}c)+\cdots \\
  & = & \phi\circ J_C(\phi^{-1}a,\phi^{-1}b,\phi^{-1}c) = 0,
\end{eqnarray*}     
where we introduced the {\sl Jacobian tensor} $J_C\in(V^*)^{\otimes 3}\otimes V=T^{1,3}V$ for an element $C\in T^{1,2}V$ by
\[
 J_C(a,b,c):=C(C(a,b),c)+C(C(b,c),a)+C(C(c,a),b).
\]
For $C\in T^{1,2}V$ skew symmetric we have that $J_C\in T^{1,3}V$ is skew symmetric as well and we look at $J_C$ therefore as an element in $\Lambda^3V^*\otimes V$. 

Let us choose a basis ${\cal A}=(e_1,\ldots,e_n)$ of $V$, e.g., the canonical basis in the case $V=K^n$. Then we can describe $C$ by its so called {\sl structure constants}
\[
 C_{ij}^k=\iprod{\lambda^k}{C(e_i,e_j)},\qquad\mbox{for $1\leq i,j,k\leq n$},
\]
where $(\lambda^1,\ldots,\lambda^n)$ denotes the dual basis of ${\cal A}$ and $\iprod{\cdot}{\cdot}$ the natural pairing between $V^*$ and $V$. Then the vanishing of $J_C\in\Lambda^3V^*\otimes V$ gives altogether
\[
 {n\choose 3}n=\frac{1}{6}(n-2)(n-1)n^2
\]
quadratic equations on the space $W=\Lambda^2V^*\otimes V\cong K^N$ of dimension
\[
 N={n\choose 2}n=\frac{1}{2}(n-1)n^2.
\] 
These equations will be in general not independent in the sense that the dimension of $X$ is given by
\[
 {n\choose 2}n-{n\choose 3}n,
\]
since this number will be negative for large $n$. However, in the special case $n=3$ where $\dim W=9$, we will see that the {\sl Lie variety} $X$ is in fact $6$-dimensional. (It would be interesting to know the points of $X$ where the differentials of the quadratic polynomials are linearly independent or all singularities of $X$ since these sets are $G$ invariant and therefore a union of orbits of $G$.)

 \item[(ii)] Then we have to clarify the orbit structure of $X$, i.e., to find all orbits in $X$ and moreover, if possible, to answer the question which orbits lie in the closure of a given orbit (with respect to the usual topology if $K=\R$ or $K=\C$ or with respect to the Zarisky topology in all cases of $K$) or even more to discuss the question of the topology of the orbit space $X/G$ (with respect to the quotient topology).

\end{description}
\bigskip

\noindent
{\bf 1.4.} An important step is now to decompose the representation $W=\Lambda^2V^*\otimes V$ into irreducible components, if this is possible. So let us try to find a direct sum of $W$ into $G$ invariant subspaces $W_1$ and $W_2$, i.e.,
\[
 W=W_1\oplus W_2,
\]
in a way such that $W_2$ is not decomposable furthermore and may be even $W_1$. We will see that this is possible in case of dimension $n=3$ and $K=\R$ or $K=\C$ (and in case $n<3$ anyway). 

For this we consider the following natural trace map,
\[
 \tr\colon\Lambda^2V^*\otimes V\to V^*,
\]
given by
\[
 \iprod{\tr(C)}{a}:=\iprod{\lambda^i}{C(e_i,a)},
\]
where $(e_i)$ is again a basis of $V$ and $(\lambda^i)$ its dual of $V^*$ (and we use the Einstein summation convention where we sum up over all $i$ in a formula, if $i$ appears twice there, once upstairs and once downstairs). This is independent of the chosen basis, since it is induced in a coordinate free way from the linear mapping
\begin{eqnarray*}
 F\colon V^*\otimes V^*\otimes V & \to & V^* \\
 \alpha\otimes\beta\otimes a & \mapsto & \iprod{\alpha}{a}\beta - \iprod{\beta}{a}\alpha,
\end{eqnarray*}
(i.e., $F=\tr_{(1,3)}-\tr_{(2,3)}$) by the universal property via the natural projection $\pi\colon T^{1,2}V\to\Lambda^2V^*\otimes V$: $-2\tr\circ\pi=F$. If we denote for every $a\in V$ as usual $\ad_a\colon V\to V$,
\[
 \ad_a(b)=C(a,b),
\]
then we see that
\[
 \iprod{\tr(C)}{a}=\iprod{\lambda^i}{C(e_i,a)}=-\iprod{\lambda^i}{\ad_a(e_i)}=-\tr(\ad_a),
\]
now meaning the usual trace of an endomorphism of $V$, and therefore $\tr(C)=0$ if and only if $\tr(\ad_a)=0$, for all $a\in V$, i.e., if $C$ is {\sl unimodular}, as one says. 

This homomorphism is $G$ equivariant, i.e.,
\[
 \tr(\phi.C)=\phi.\tr(C),
\]
for all $C\in W$ and $\phi\in G$, and it has, up to a factor $(1-n)$, a $G$ equivariant section $j\colon V^*\to\Lambda^2V^*\otimes V$, i.e.,
\begin{equation} \label{eqn:tr}
 \tr\circ j=(1-n)\id,
\end{equation}
which we introduce now. Consider
\[
 j\colon V^*\to\Lambda^2V^*\otimes V\cong\Alt_2(V,V)
\]
(denoting by $\Alt_2(V,U)$ the alternating bilinear forms on $V$ with values in $U$) given by
\[
 j(\alpha)(a,b)=\iprod{\alpha}{a}b-\iprod{\alpha}{b}a.
\]
Then, in fact,
\begin{eqnarray*}
 \iprod{\tr\circ j(\alpha)}{a} & = & \iprod{\lambda^i}{j(\alpha)(e_i,a)} = \iprod{\lambda^i}{\iprod{\alpha}{e_i}a-\iprod{\alpha}{a}e_i} \\
 & = & \iprod{\alpha}{\underbrace{\lambda^i(a)e_i}_{=a}} - \iprod{\alpha}{a}\underbrace{\iprod{\lambda^i}{e_i}}_{=n} = (1-n)\iprod{\alpha}{a}
\end{eqnarray*}
for all $\alpha\in V^*$ and $a\in V$, i.e.,
\[
 \tr\circ j=(1-n)\id.
\]
So let us assume that $n>1$ and the characteristic of $K$ is zero here, $\char(K)=0$. Then $(1-n)^{-1}j$ is in fact a section of $\tr$. 

The case $n=1$ is trivial, i.e., there exists only the trivial structure (since $W=\Lambda^2V^*\otimes V=(0)$, if you want) and in case $n=2$, where $\dim W=\dim V^*=2$, we see that $\tr$ (and $j$) is already an isomorphism between $G$ representations (or, as we also say, $G$ modules). But $V^*$ has only two $G$ orbits, namely $\{0\}$ and $V^*\setminus\{0\}$. (In fact, if $\alpha,\beta\in V^*\setminus\{0\}$, we can choose $a,b\in V$ such that
\[
 V=\ker(\alpha)+\langle a\rangle=\ker(\beta)+\langle b\rangle
\]
and, by normalizing $a,b\in V$ if necessary, 
\[
 \alpha(a)=\beta(b)=1.
\]
Now we can find a $\phi\in G$ with $\phi(\ker\alpha)=\ker\beta$ and $\phi(a)=b$. Then
\[
 \phi^*\beta|\ker\alpha=\beta\circ\phi|\ker\alpha=\beta|\ker\beta=0
\]
and
\[
 \phi^*\beta(a)=\beta\circ\phi(a)=\beta(b)=\alpha(a)
\]
and therefore $\phi^*\beta=\alpha$.) Therefore, there exist only two Lie structures in this dimension. The {\sl abelian} structure, i.e., $C=0$, is the only unimodular here, and the {\sl solvable} structure, since
\[
 V^{(2)}:=[[V,V],V]=0.
\]
\bigskip

\noindent
{\bf 1.5.} In this way we have decomposed $W=\Lambda^2V^*\otimes V$ into the direct sum of invariant subspaces $W_1:=\ker(\tr)$ and $W_2=j(V^*)$, where $W_2$ is irreducible, i.e., the only invariant subspaces of $W_2$ are the trivial ones. This is the case since $W_2$ is just a copy of $V^*$ which has only two orbits as we have seen. An invariant subspace $U\ne(0)$ must therefore be the whole space. Thus we have to investigate $W_1$ now. 

For this we introduce a similar extension construction as we considered it for the module $V^*=\Lambda^1V^*\otimes S^0V$ to the module $W=\Lambda^2V^*\otimes S^1V$, namely
\[
 p\colon\Lambda^2V^*\otimes V\to\Lambda^3V^*\otimes S^2V \cong \Alt_3(V,\Sym_2(V^*)),
\]
\begin{eqnarray*}
 \lefteqn{p(C)(a,b,c)(\alpha,\beta) :=} \\
  & & \iprod{\alpha}{C(a,b)}\iprod{\beta}{c} + \iprod{\beta}{C(a,b)}\iprod{\alpha}{c} \\
  & & \iprod{\alpha}{C(b,c)}\iprod{\beta}{a} + \iprod{\beta}{C(b,c)}\iprod{\alpha}{a} \\
  & & \iprod{\alpha}{C(c,a)}\iprod{\beta}{b} + \iprod{\beta}{C(c,a)}\iprod{\alpha}{b}. 
\end{eqnarray*}
$p$ is in fact skew symmetric in $a,b,c$ and symmetric in $\alpha,\beta$. We call it $p$ for projection (and not like $j$ for inclusion) because in our case $n=3$ we will see that $p$ is surjective and restricted to $W_1=\ker(\tr)$ even an isomorphism.

But observe first
\begin{eqnarray*}
 \lefteqn{p\circ j(\alpha)(a,b,c)(\beta,\gamma)=\iprod{\beta}{\alpha(a)b-\alpha(b)a}\gamma(c)+\cdots} \\
 & = & \alpha(a)\beta(b)\gamma(c)-\alpha(b)\beta(a)\gamma(c) \\
 &   & +\cdots \\
 & = & 0,
\end{eqnarray*}
for all $\alpha,\beta,\gamma,a,b,c$, i.e.,
\begin{equation} \label{eqn:pj}
 p\circ j=0.  
\end{equation}
\bigskip

\noindent
{\bf 1.6.} Now we consider, in analogy to $\tr\colon W\to V^*$, the trace map
\[
 \Tr\colon\Lambda^3V^*\otimes S^2V\to\Lambda^2V^*\otimes V\cong\Alt_2(V,V),
\]
\[
 \iprod{\alpha}{\tr(M)(a,b)}:=M(a,b,e_i)(\lambda^i,\alpha)
\]
for all $a,b\in V$ and $\alpha\in V^*$ (and $(e_i)$ a basis of $V$ and $(\lambda^i)$ its dual basis of $V^*$). Then we compute for all $a\in V$
\begin{eqnarray*}
 \iprod{\tr\circ\Tr(M)}{a} & = & \iprod{\lambda^i}{\Tr(M)(e_i,a)} \\
=M(e_i,a,e_j)(\lambda^j,\lambda^i) & = & -M(e_j,a,e_i)(\lambda^i,\lambda^j)
\end{eqnarray*}
(because of the skew symmetry in $(e_i,a,e_j)$ and symmetry in $(\lambda^i,\lambda^j)$)
\begin{eqnarray*}
 & = & -M(e_i,a,e_j)(\lambda^j,\lambda^i)\qquad\mbox{(by permutation of the summation indices)} \\
 & = & -\iprod{\tr\circ\Tr(M)}{a},
\end{eqnarray*}
i.e.,
\begin{equation} \label{eqn:trace}
 \tr\circ\Tr=0.
\end{equation}
\bigskip

\noindent
{\bf 1.7.} Now let us assume that $\dim V=3$.

\begin{lemma} Then we have:
 \begin{equation} \label{eqn:p_circ_Tr}
 p\circ\Tr=2\id.
 \end{equation}
\end{lemma}

\noindent
\begin{proof} Let $\lambda^3\in V^*\setminus\{0\}$ be arbitrary. It is enough to show for an arbitrary basis $(e_1,e_2,e_3)$ of $V$ that
\[
 p\circ\Tr(M)(e_1,e_2,e_3)(\lambda^3,\lambda^3)=2M(e_1,e_2,e_3)(\lambda^3,\lambda^3)
\]
is true. For this implies 
\[
 p\circ\Tr(M)(a,b,c)(\lambda^3,\lambda^3)=2M(a,b,c)(\lambda^3,\lambda^3),
\]
for all $a,b,c\in V$ and since $\lambda^3\in V^*\setminus\{0\}$ was arbitrary we conclude by the symmetry of $M(a,b,c)$ that
\[
 p\circ\Tr(M)=2M,
\]
for all $M\in\Lambda^3V^*\otimes S^2V$.

Now we complement $\lambda^3$ to a basis $(\lambda^1,\lambda^2,\lambda^3)$ of $V^*$ and let $(e_1,e_2,e_3)$ be the dual to it. Then we compute indeed
\begin{eqnarray*}
 \lefteqn{\frac{1}{2}p\circ\Tr(M)(e_1,e_2,e_3)(\lambda^3,\lambda^3)} \\
 & = & \iprod{\lambda^3}{\tr(M)(e_1,e_2)}\lambda^3(e_3) + \iprod{\lambda^3}{\tr(M)(e_2,e_3)}\lambda^3(e_1) \\
 &  & + \iprod{\lambda^3}{\tr(M)(e_3,e_1)}\lambda^3(e_2) \\
 & = & \iprod{\lambda^3}{\tr(M)(e_1,e_2)} = M(e_1,e_2,e_i)(\lambda^i,\lambda^3) \\
 & = & M(e_1,e_2,e_3)(\lambda^3,\lambda^3).
\end{eqnarray*}
\end{proof}
\bigskip

\noindent
In particular, we see that $p$ is in fact surjective and $\Tr$ is injective. Since we already know that $\im(\Tr)\sub W_1$ and
\[
 \dim(\Lambda^3V^*\otimes S^2V)=6=\dim W_1,
\]
$\Tr\colon\Lambda^3V^*\otimes S^2V\to W_1$ must be even an isomorphism (of $G$ modules since $\Tr$ is equivariant as well). The same is true for $p|W_1\colon W_1\to\Lambda^3V^*\otimes S^2V$.

\begin{corollary} \label{cor:decomp}
Let $\dim V=3$. Then:
\[
 \frac{1}{2}\Tr\circ p - \frac{1}{2}j\circ\tr=\id.
\]
\end{corollary}

\noindent
\begin{proof} Using the equations (\ref{eqn:tr})-(\ref{eqn:pj}) we get
\[
 (\frac{1}{2}\Tr\circ p - \frac{1}{2}j\circ\tr)\circ j=-\frac{1}{2}j\circ\tr\circ j=j
\]
and therefore
\[
 (\frac{1}{2}\Tr\circ p - \frac{1}{2}j\circ\tr)|W_2=\id_{W_2}
\]
since $W_2=\im(j)=\ker(p)$. On the other hand, now using equations (\ref{eqn:trace})-(\ref{eqn:p_circ_Tr}), we have
\[
 (\frac{1}{2}\Tr\circ p - \frac{1}{2}j\circ\tr)\circ\Tr=\frac{1}{2}\Tr\circ p\circ\Tr=\Tr,
\]
and therefore also
\[
 (\frac{1}{2}\Tr\circ p - \frac{1}{2}j\circ\tr)|W_1=\id_{W_1}
\]
since $W_1=\im(\Tr)=\ker(\tr)$. From
\[
 W=W_1\oplus W_2
\]
we get the assertion.
\end{proof}
\bigskip

\noindent
Thus we can summarize our investigations in the following two exact $G$-sequences of $G$-modules,
\[
 0\longrightarrow W^\prime\stackrel{\Tr}{\longrightarrow}W\stackrel{\tr}{\longrightarrow}V^*\longrightarrow 0,
\]
\[
 0\longleftarrow W^\prime\stackrel{p}{\longleftarrow}W\stackrel{j}{\longleftarrow}V^*\longleftarrow 0,
\]
together with the equations
\begin{eqnarray*}
 \tr\circ j & = & -2\id \\
 p\circ\Tr & = & 2\id \\
 \Tr\circ p-j\circ\tr & = & 2\id.
\end{eqnarray*}
We have therefore decomposed our $G$-module $W=\Lambda^2V^*\otimes V$ into the direct sum of the invariant subspaces $W_1=\Tr(\Lambda^3V^*\otimes S^2V)$ and $W_2=j(V^*)$,
\[
 W=W_1\oplus W_2.
\]
For $C\in W$ we call $\frac{1}{2}\Tr\circ p(C)\in W_1$ the {\sl unimodular (or tracefree) part of $C$} and $-\frac{1}{2}j\circ\tr(C)\in W_2$ the {\sl trace part of $C$}.

$W_2\cong V^*$ is an irreducible $G$-module as we have seen, and when we have studied the orbit structure of $W_1\cong\Lambda^3V^*\otimes S^2V$ we will see that $W_1$ is irreducible as well. From the orbit structures of $W_1$ and $W_2$ (and the isotropy groups of $W_1$) we will get the orbit structure of $W$, more precisely the orbit structure of the Lie variety $X\sub W$.
\bigskip

\noindent
{\bf 1.8.} Before doing this let us consider now the Jacobian identity of a Lie structure. Remember that for an arbitrary element $C\in\Lambda^2V^*\otimes V$ we defined its Jacobian tensor $J_C\in\Lambda^3V^*\otimes V\cong\Alt_3(V,V)$ by
\[
 J_C(a,b,c):=C(C(a,b),c)+\,\mbox{c.p.}
\]
(c.p.\ abbreviating cyclic permutations). Thus our Lie variety is given by
\[
 X=\{C\in\Lambda^2V^*\otimes V:\; J_C=0\}.
\]
Now let $M\in\Lambda^3V^*\otimes S^2V$ and $\nu\in V^*$. Then we consider
\[
 C=\Tr(M)+j(\nu),
\]
where we know from \ref{cor:decomp} that every $C\in\Lambda^2V^*\otimes V$ is of this form, $M$ and $\nu$ uniquely determined by
\[
 M=\frac{1}{2}p(C)\qquad\nu=-\frac{1}{2}\tr(C).
\]
For an arbitrary $\alpha\in V^*$ we denote by $i_\alpha M\in\Lambda^3V^*\otimes V$ the {\sl contraction of $M$ in the $4$. (or $5$.) argument}, i.e.,
\[
 i_\alpha M(a,b,c)=M(a,b,c)(\alpha,\cdot),
\]
or, if we choose a basis $(e_i)$ of $V$ with dual $(\lambda^i)$:
\[
 i_\alpha M(a,b,c)=M(a,b,c)(\alpha,\lambda^1)e_1+\,\cdots+M(a,b,c)(\alpha,\lambda^3)e_3.
\]

\begin{theorem} Let $\dim V=3$. Then for every $M\in\Lambda^3V^*\otimes S^2V$ and $\nu\in V^*$ we have for $C:=\Tr(M)+j(\nu)$:
\[
 J_C=2i_\nu M.
\]
\end{theorem}

\noindent
\begin{proof} Let $(e_1,\ldots,e_3)$ be a basis of $V$ and $(\lambda^1,\ldots,\lambda^3)$ be its dual basis of $V^*$. By definition for $C=\Tr(M)+j(\nu)$ we have
\begin{eqnarray*}
 \lefteqn{J_C(e_1,\ldots,e_3)=C(C(e_1,e_2),e_3)+\,\mbox{c.p.}} \\
 & = & C(\Tr M(e_1,e_2),e_3) + C(j(\nu)(e_1,e_2),e_3) +\,\mbox{c.p.}
\end{eqnarray*}
Now the 1.\ summand (together with its cyclic permutations) vanishes since
\begin{eqnarray*}
 \lefteqn{C(\Tr M(e_1,e_2),e_3)} \\
 & = & C(M(e_1,e_2,e_i)(\lambda^i,\lambda^1)e_1)+\ldots+M(e_1,e_2,e_i)(\lambda^i,\lambda^3)e_3,e_3) \\
 & = & M(e_1,e_2,e_3)(\lambda^3,\lambda^1)C(e_1,e_3) + M(e_1,e_2,e_3)(\lambda^3,\lambda^2)C(e_2,e_3),
\end{eqnarray*} 
because $C(e_3,e_3)=0$ and $M(e_1,e_2,e_i)=0$ for $i=1,2$. Therefore we get with $M_{ijk}:=M(e_i,e_j,e_k)$:
\begin{eqnarray*}
 \lefteqn{C(\Tr M(e_1,e_2),e_3)+\,\mbox{c.p.}} \\
 & = & M_{123}(\lambda^3,\lambda^1)C(e_1,e_3)+M_{123}(\lambda^3,\lambda^2)C(e_2,e_3) \\
 &   & +M_{231}(\lambda^1,\lambda^2)C(e_2,e_1)+M_{231}(\lambda^1,\lambda^3)C(e_3,e_1) \\
 & & +M_{312}(\lambda^2,\lambda^3)C(e_3,e_2)+M_{312}(\lambda^2,\lambda^1)C(e_1,e_2) \\
 & = & 0,
\end{eqnarray*}
because of the skew symmetry of $C$ and the symmetry of $M$ (in $V^*\times V^*$). (From this we see already that $J_C=0$ for $\nu=0$, i.e., $W_1\sub X$.)

For the 2.\ summand we have
\[
 C(j(\nu)(e_1,e_2),e_3)=\Tr M(j(\nu)(e_1,e_2),e_3)+j(\nu)(j(\nu)(e_1,e_2),e_3),
\]
from which the second summand (with its permutations) vanishes (in any dimension) since
\begin{eqnarray*}
 \lefteqn{j(\nu)(j(\nu)(e_1,e_2),e_3)+\,\mbox{c.p.} = j(\nu)(\nu(e_1)e_2-\nu(e_2)e_1, e_3)+\,\mbox{c.p}} \\
 & = & \nu(e_1)\nu(e_2)e_3-\nu(e_1)\nu(e_3)e_2-\nu(e_2)\nu(e_1)e_3+\nu(e_2)\nu(e_3)e_1+\,\mbox{c.p} \\
 & = &\nu(e_2)\nu(e_3)e_1-\nu(e_3)\nu(e_1)e_2+\,\mbox{c.p.} \\
 & = & 0.
\end{eqnarray*}
Finally we compute for the first summand
\begin{eqnarray*}
 \lefteqn{\Tr M(j(\nu)(e_1,e_2),e_3)+\,\mbox{c.p.}} \\
 & = & \nu(e_1)\Tr M(e_2,e_3)-\nu(e_2)\Tr M(e_1,e_3) +\,\mbox{c.p.} \\
 & = & \nu(e_1)\Tr M(e_2,e_3)+\nu(e_2)\Tr M(e_3,e_1) +\,\mbox{c.p.} \\
 & = & 2(\nu(e_1)\Tr M(e_2,e_3)+\ldots+\nu(e_3)\Tr M(e_1,e_2)) \\
 & = & 2 M_{123}(\nu(e_1)\lambda^1+\ldots+\nu(e_3)\lambda^3,\cdot) \\
 & = & 2 M_{123}(\nu,\cdot) = 2 i_\nu M(e_1,e_2,e_3),
\end{eqnarray*}
and therefore in fact
\[
 J_C=2i_\nu M.
\]
\end{proof}
\bigskip

\noindent
Although the algebraic variety $X\sub W\cong K^9$ is described by three quadratic equations (namely the components of $J_C$ w.r.t.\ a choice of basis of $V$), it has w.r.t.\ the projection $p\colon W\to W^\prime$ some remarkable linear structures. First of all the zero section $W_1=\Tr(W^\prime)$ with $W^\prime=\Lambda^3V^*\otimes S^2V$ of $p$ is a component of $W$ in the algebraic geometric sense (since $\dim W_1=\dim X=6$, as we will see). On the other hand the fiber of any point $M\in W^\prime$ is (up to a translation with $M$)  a linear subspace of $W$, more precisely of $V^*$ (embedded in $W$ via $j$), namely the nullspace of the symmetric bilinear form 
\[
 M\in\Lambda^3V^*\otimes S^2V\cong\Sym_2(V^*,\Lambda^3V^*).
\]
For instance the subspace $W_2=j(V^*)$ belongs to $X$ since it is the nullspace of $M=0$, of course.

We will see soon that $W^\prime$ consists in the real case of exactly $6$ orbits
(and in the complex case of $4$ orbits), $2$ of them being open in $W^\prime$ (only $1$ in the complex case), corresponding to those forms which are non degenerate. The complement consists of a $5$-dimensional cone (as the zero set of a determinant) and only over their orbits, besides the origin, the fibers of $p$ consists of non trivial (affine) subspaces of $W$ with translation vector space in $j(V)$.
\bigskip

\noindent
{\bf 1.9.} Let us now investigate the orbit structure of the $G$-module $W^\prime:=\Lambda^3V^*\otimes S^2V$. For that we choose a basis $(e_i)$ of $V$ and let as usual be $(\lambda^i)$ dual to it. An element $M\in W^\prime$ is then described by a matrix $A=(a^{ij})$ via
\[
 a^{ij}:=M(e_1,e_2,e_3)(\lambda^i,\lambda^j).
\]
Let further be $\phi\in G=\Aut(V)$. Then $\phi$ is described by the matrix $g=(g^i_j)$ via
\[
 g^i_j:=\iprod{\lambda^i}{\phi(e_j)}.
\]
The action of $G$ on $W^\prime$, namely
\[
 (\phi.M)(a,b,c)(\alpha,\beta)=M(\phi^{-1}a,\phi^{-1}b,\phi^{-1}c)(\phi^{-1}\alpha,\phi^{-1}\beta),
\] 
translates into the following action of $\GL_n(K)$ on
\[
 \Sym_2(K):=\{A\in\Mat_2(K):\; A^T=A\}:
\]
If $B=(b^{ij})$ describes $\phi.M$ w.r.t.\ $(e_i)$, then
\begin{eqnarray*}
 b^{ij} & = & (\phi.M)(e_1,e_2,e_3)(\lambda^i,\lambda^j) \\
 & = & M(\phi^{-1}e_1,\phi^{-1}e_2,\phi^{-1}e_3)(\phi^{-1}\lambda^i,\phi^{-1}\lambda^j) \\
 & = & \det(\phi^{-1})M(e_1,e_2,e_3)(\phi^{-1}\lambda^i,\phi^{-1}\lambda^j).
\end{eqnarray*}
From the definition of the matrix $g$ we get immediately that
\[
 \phi(e_i)=\iprod{\lambda^k}{\phi(e_i)}e_k=g^k_ie_k,
\]
and from the action of $G$ on $V^*$, given by
\[
 \iprod{\phi.\lambda}{a}=\iprod{\lambda}{\phi^{-1}a},
\]
we see that
\[
 \iprod{\phi^{-1}\lambda^i}{e_j}=\iprod{\lambda^i}{\phi.e_j}=\iprod{\lambda^i}{g^k_je_k}=g^i_j,
\]
i.e.,
\[
 \phi^{-1}.\lambda^i=\iprod{\phi^{-1}.\lambda^i}{e_k}\lambda^k=g^i_k\lambda^k.
\]
From this we deduce that
\begin{eqnarray*}
 b^{ij} & = & \det(g^{-1})g^i_kg^j_lM_{123}(\lambda^k,\lambda^l) \\
 & = &\det(g^{-1})g^i_ka^{kl}g^j_l,
\end{eqnarray*}
or, in matrix description,
\[
 B=\det(g^{-1})gAg^T.
\]
\bigskip

\noindent
Let us assume now $K=\R$ and set $U=\Sym_2(\R^3)$. By Sylvester's classical theorem the orbits of $U$ under the natural action
\begin{equation} \label{eqn:Sylvester}
 g.A=gAg^T
\end{equation}
of $\GL_3\R$ ($g\in\GL_3\R$, $A\in\Sym_2(\R^3)$) are classified by the rank $r=\rk(A)\in\N_0$, and its signature $s=\sgn(A)\in\{-r+2k\in\Z:$ $k=0,\ldots r\}$. Therefore there exist $10$ orbits described by triples $(k,l,m)$, where $k$ is the dimension of a maximal subspace of $\R^3$, on which $A$ is positive definite, $l$ the dimension of a maximal subspace where $A$ is negative definite and $m$ ist the dimension of the null space of $A$:
\[
 \begin{array}{cccccc}
 (1,1,1)   & (1,1,-1)  & (1,1,0)   & (1,-1,0) & (1,0,0)  & (0,0,0) \\
 (-1,-1,-1)& (-1,-1,1) & (-1,-1,0) &          & (-1,0,0) &
 \end{array}
\]
or in terms of $(r,s)$, using $r=k+l$ and $s=k-l$,
\[
 \begin{array}{cccccc}
 (3,3)  & (3,1)  & (2,2)  & (2,0) & (1,1)  & (0,0) \\
 (3,-3) & (3,-1) & (2,-2) &       & (1,-1) &
 \end{array}
\]
Any of the orbits of this natural action on $\Sym_2(\R^3)$ is now contained in an orbit of the slightly more complicated action given by
\begin{equation} \label{eqn:W_1-action}
 g.A=\frac{1}{\det(g)}gAg^T
\end{equation}
of $G=\GL_3\R$ on $U=\Sym_2(\R^3)$, i.e., the orbits under the action of (\ref{eqn:W_1-action}) are unions of orbits under the action of (\ref{eqn:Sylvester}). Indeed, if 
\[
 B=gAg^T
\]
for some $g\in G$ (and given $A\in U$), consider $\tilde g:=\frac{1}{\det(g)}g$ and compute 
\[
 \det(\tilde g)=\det(g)^{-3}\det(g)=\det(g)^{-2}
\]
and therefore
\[
 \frac{1}{\det(\tilde g)}\tilde gA\tilde g^T=\det(g)^2\frac{1}{\det(g)}gAg^T\frac{1}{\det(g)}=gAg^T.
\]
On the other hand two orbits $B_1$ and $B_2$ under the action of (\ref{eqn:Sylvester}) can belong to the same orbit under the action of (\ref{eqn:W_1-action}) only if $B_2=\pm B_1$, since the multiplication with $\det(g)^{-1}$ cannot change the rank of $A$ and from the signature at most its sign. The {\sl absolute value of the signature} is therefore an invariant under the action of (\ref{eqn:W_1-action}). In fact, if $A$ lies in an orbit, then $-A$ lies in the same orbit as can be seen by acting with the element $g=-\1$:
\[
 g.A=\frac{1}{\det(g)}gAg^T=-A.
\]

\begin{proposition}
On the $G$ module $W^\prime=\Lambda^3V^*\otimes S^2V$ there exist exactly $6$ orbits. They are classified by their rank and the absolute value of its signature.
\end{proposition}

\noindent
We leave it to the reader to give a coordinate free definition of the rank and the absolute value for an element
\[
 M\in\Lambda^3V^*\otimes S^2V\cong\Sym_2(V^*,\Lambda^3V^*).
\]
For the case $\R^3$ let us fix for every orbit the following natural representatives of the orbits:
\[
\begin{array}{ll}
 \Bia(I):     &  A=\diag(0,0,0)  \\
 \Bia(II):    &  A=\diag(1,0,0)  \\
 \Bia(VI_0):  &  A=\diag(1,-1,0) \\
 \Bia(VII_0): &  A=\diag(1,1,0)  \\
 \Bia(VIII):  &  A=\diag(1,1,-1) \\
 \Bia(IX):    &  A=\diag(1,1,1)
\end{array}
\]
We want to mention here that $\Bia(I)$ is the {\sl abelian} Lie algebra, $\Bia(II)$ is the (unique) {\sl nilpotent} (and non abelian) Lie algebra, $\Bia(VI_0)$ and $\Bia(VII_0)$ are {\sl solvable} (and non nilpotent, of course), and $\Bia (VIII)$ and $\Bia(IX)$ are the simple Lie algebras over the reals in dimension $3$.

In a similar way one sees that the orbits in the complex case are just the same as for the natural action and are classified just by its rank. So there are altogether $4$ orbits given by the representives $\Bia(I)$ (the abelian), $\Bia(II)$ (the nilpotent), $\Bia(VII_0)$ (the solvable) and $\Bia(IX)$ (the simple).

Finally we want to refer to J.\ Milnors beautiful, and more elementary, proof of the classification of the unimodular Lie algebras in dimension $3$ (cf.\ \cite{Milnor}).

\section{Classification of Bianchi B classes}

%Abschnitt 1: Bianchi-Klassifikation der 3-dimensionalen Liealgebren, Klasse A%

{\bf 2.1.} We want to start now with the computation of the isotropy groups in the $6$ representatives ($4$ in the complex case) of the $3$-dimensional unimodular Lie algebras. This will give us at first the dimensions of the $6$ orbits. These carry in fact a manifold structure via the natural orbit maps $\rho\colon G\to W_1$,
\[
 \rho(g)=g.A_0
\]
($A_0$ one of the representatives from the last paragraph), since these induce a smooth and injective immersion $\bar\rho\colon G/H\to W_1$ with image $B=G(A_0)$, the orbit of $A_0$, and $G/H$ has a natural manifold structure of dimension
\[
 \dim(G/H)=\dim G-\dim H.
\]
Next these isotropy groups 
\[
 H=G_{A_0}=\{g\in G:\; g.A_0=A_0\}
\]
will give us the {\sl automorphism groups} of the $6$ (resp.\ $4$) unimodular Bianchi algebras since obviously $H\sub G$ consists exactly of those vector space isomorphisms respecting the Lie structure.

Finally the computation will also be the key for the classification of the non unimodular Lie groups (so called $B$-class in the terminology of Bianchi).

\begin{proposition}
Let $n\in\N$ and $(k,l,m)\in\N_0^3$ be such that $k+l+m=n$ and let
\[
 A_0=\diag(\underbrace{1,\ldots,1}_{\mbox{$k$ times}},\underbrace{-1,\ldots,-1}_{\mbox{$l$ times}},\underbrace{0,\ldots,0}_{\mbox{$m$ times}}).
\]
Let furthermore act $G=\GL_n\R$ on $W^\prime=\Sym_2(\R^n)$ via
\[
 g.A=\det(g^{-1})gAg^T.
\]
Then the isotropy group $H=G_{A_0}$ consists of all invertible matrices of the form
\begin{equation} \label{eqn:matrix}
 g=\left(\begin{array}{cc} g_1 & h \\ 0 & g_2 \end{array}\right)
\end{equation}
with $g_1\in\GL_{k+l}\R$, $g_2\in\GL_m\R$, and $h\in\Mat(k+l,m;\R)$ satisfying
\begin{equation} \label{eqn:determinant}
 g_1\1_{k,l}g_1^T=\det(g_1)\det(g_2)\1_{k,l}.
\end{equation}
\end{proposition}
Here we use the notation
\[
 \1_{k,l}=\diag(\underbrace{1\ldots,1}_{\mbox{$k$ times}},\underbrace{-1,\ldots,-1}_{\mbox{$l$ times}}).
\]

\noindent
\begin{proof} The condition $g.A_0=A_0$ reads in block form
\[
 \det(g^{-1})\left(\begin{array}{cc} g_1 & h_1 \\ h_2 & g_2 \end{array}\right)\left(\begin{array}{cc} \1_{k,l} & 0 \\ 0 & 0 \end{array}\right)\left(\begin{array}{cc}g_1^T & h_2^T \\ h_1^T & g_2^T \end{array}\right)=\left(\begin{array}{cc}\1_{k,l} & 0 \\ 0 & 0 \end{array}\right),
\]
where we have set
\[
 g=\left(\begin{array}{cc} g_1 & h_1 \\ h_2 & g_2 \end{array}\right).
\]
This comes down to
\[
 \det(g)\left(\begin{array}{cc} \1_{k,l} & 0 \\ 0 & 0 \end{array}\right) = \left(\begin{array}{cc} g_1\1_{k,l}g_1^T & g_1\1_{k,l}h_2^T \\ h_2\1_{k,l}g_1^T & h_2\1_{k,l}h_2^T \end{array}\right).
\]
Looking into the left top corner we conclude at first that
\[
 \det(g)\1_{k,l}=g_1\1_{k,l}g_1^T.
\]
Taking determinants we conclude that
\[
\det(g)^{k+l}(-1)^l=\det(g_1)^2(-1)^l,
\]
from which we see that $\det(g_1)\ne 0$ since $\det(g)\ne 0$. Turning to the right top corner we can conclude the following:
\[
 g_1\1_{k,l}h_2^T=0\;\Rightarrow\;h_2=0
\]
(by multiplying with $\1_{k,l}g_1^{-1}$ from left). The conditions from the bottom row are then automatically fulfilled.

Next we see that
\[
 \det(g)=\det(g_1)\det(g_2),
\]
which implies $\det(g_2)\ne 0$ as well and the condition of the left top corner becomes
\[
 g_1\1_{k,l}g_1^T=\det(g_1)\det(g_2)\1_{k,l}.
\]
This shows that the elements of the isotropy groups fulfill the conditions (\ref{eqn:matrix}) and (\ref{eqn:determinant}).

Vice versa, if $g\in\GL_n\R$ satisfy (\ref{eqn:matrix}) and (\ref{eqn:determinant}), then our computation shows that $g$ lies in fact in the isotropy group of $A_0$.
\end{proof}
\bigskip

\noindent
Observe next that the representation of $G=\Aut(V)$ on 
\[
 W^\prime=\Lambda^3V^*\otimes S^2V\cong\Sym_2(V^*,\Lambda^3V^*)
\]
induces a representation of the isotropy group $H=G_M$, for some element $M\in W^\prime$, on the {\sl null space} 
\[
 U=\Deg(M)=\{\alpha\in V^*:\; M(\alpha,\beta)=0,\,\mbox{for all $\beta\in V^*$}\}
\]
by the homomorphism $\pi\colon G\to\Aut(U)$,
\[
 \pi(h)=h^*|U.
\]
In fact, $U\sub V^*$ is invariant under $h^*$ since for $\alpha\in U$, i.e., $M(\alpha,\beta)=0$ for all $\beta\in V^*$, we compute
\[
 M(h^*\alpha,\beta)=M(h^{-1}.\alpha,h^{-1}.(h.\beta)),
\]
since $h^{-1}.\alpha=\alpha\circ h=h^*(\alpha)$, and therefore
\[
 M(h^*\alpha,\beta)=h^{-1}.M(\alpha,h.\beta)=0,
\]
by the equivariance of $M$. This is true for any $\beta\in V^*$, i.e., $h^*\alpha\in U$ implying $\pi(h)\in\Aut(U)$.

W.r.t.\ the Bianchi representatives of the unimodular structures on $V=\R^3$ this comes down simply to the projection of $g\in H$ to its right lower corner,
\[
 \pi\colon H\to\GL_m\R,\;\left(\begin{array}{cc} g_1 & h \\ 0 & g_2 \end{array}\right)\mapsto g_2.
\]
Next observe, by taking determinants of equation (\ref{eqn:determinant}), we get the equation
\[
 \det(g_1)^2(-1)^l=\det(g_1)^{k+l}\det(g_2)^{k+l}(-1)^l,
\]
i.e.,
\begin{equation} \label{eqn:Det}
 \det(g_1)^{2-(k+l)}=\det(g_2)^{k+l}.
\end{equation}
From this we deduce the following.

\begin{corollary} \label{cor:}
Let $H\sub\GL_n\R$ be the isotropy group of $A_0\in\Sym_n\R$ given by the triple $(k,l,m)\in\N_0^3$. Let further $\pi\colon H\to\GL_m$ be the natural projection. Then the following is true:
 \begin{description}
  \item[(a)] If $k+l\ne 2$, then $\pi$ is surjective.
  \item[(b)] If $k+l=2$, then
\[
 \pi(H)=\{g_2\in\GL_m\R:\; \det(g_2)=\pm 1\}.
\]
 \end{description}
\end{corollary}

\noindent
\begin{proof} 
(a) So let $g_2\in\GL_m(\R)$ be arbitrary. We distinguish two cases:
   \begin{enumerate} 
    \item $k+l$ is even or $\det(g_2)>0$: In this case we can take the $(2-k-l)$-th root of $\det(g_2)$,
\[
 \lambda:=\det(g_2)^{\frac{1}{2-k-l}},
\]
and we set
\[
 g_1:=\lambda\1_{k+l}.
\]
Then we compute
\begin{eqnarray*}
 \det(g_1)\det(g_2) & = & \lambda^{k+l}\det(g_2) \\ & = & \left(\det(g_2)^{\frac{1}{2-k-l}}\right)^{k+l}\cdot\left(\det(g_2)^{\frac{1}{2-k-l}}\right)^{2-k-l} \\ & = & (\det(g_2)^{\frac{1}{2-k-l}})^2=\lambda^2,
\end{eqnarray*}
and from this
\[
 g_1\1_{k,l}g_1^T=\lambda^2\1_{k,l}=\det(g_1)\det(g_2)\1_{k,l}.
\]
Therefore 
\[
 g=\left(\begin{array}{cc} g_1 & 0 \\ 0 & g_2 \end{array}\right)=:\diag(g_1,g_2)
\]
is indeed in $H$ and, of course, a preimage of $g_2$ under $\pi$.
   \item $k+l$ is even and $\det(g_2)<0$: In this case let $g_0:=\diag(-1,1,\ldots,1)$ and set $g_1:=\lambda g_0$ with
\[
 \lambda:=-(-\det(g_2))^{\frac{1}{2-k-l}}.
\]
Then we compute
\begin{eqnarray*}
 \det(g_1)\det(g_2) & = & -\lambda^{k+l}\det(g_2)=-\left(-(-\det(g_2))^{\frac{1}{2-k-l}}\right)^{k+l}\det(g_2) \\
 & = & (-\det(g_2))^{\frac{k+l}{2-k-l}}(-\det(g_2))=(-\det(g_2))^{\frac{k+l}{2-k-l} +1} \\ & = & \left((\-\det(g_2))^{\frac{1}{2-k-l}}\right)^2=\lambda^2,
\end{eqnarray*}
and from this again
\[
 g_1\1_{k,l}g_1^T=\lambda^2(g_0\1_{k,l}g_0^T)=\lambda^2\1_{k,l}=\det(g_1)\det(g_2)\1_{k,l}.
\]
Thus $g=\diag(g_1,g_2)\in H$ and $\pi$ is surjective in case (a).
   \end{enumerate}

(b) From equation (\ref{eqn:Det}) we now find
\[
 \det(g_2)^2=\det(g_1)^0=1,
\]
and therefore
\[
 \pi(H)\sub\{g_2\in\GL_m\R:\; \det(g_2)=\pm 1\}.
\]
On the other hand let $g_2\in\GL_m\R$ with $\det(g_2)=\pm 1$. We distinguish again two cases.
  \begin{enumerate} 
   \item $\det(g_2)=1$: Then we let $g_1=\1_2$, from which we get
\[
 g_1\1_{k,l}g_1^T=\1_{k,l}=\det(g_1)\det(g_2)\1_{k,l},
\]
and therefore $g=\diag(g_1,g_2)\in H$ with $\pi(g)=g_2$.
  \item $\det(g_2)=-1$: Then we set $g_1=\diag(-1,1)$ and therefore $\det(g_1)=-1$. From this again
\[
 g_1\1_{k,l}g_1^T=\1_{k,l}=\det(g_1)\det(g_2)\1_{k,l},
\]
and therefore $g_2\in\im(\pi)$.
   \end{enumerate}
 \end{proof}
\bigskip

\noindent
The complex case is quite similar. If $\rk(A_0)=r$, then $\pi\colon H\to\GL_m\C$ is surjective for $r\ne 2$. To see this set again $g_1=\lambda\1_r$ with
\[
 \lambda:=(\det(g_2))^{\frac{1}{2-r}}
\]
(for a given $g_2\in\GL_m\C$), where we take an arbitrary $(2-r)$-th root (meaning, of course, the reciprocal of an $(r-2)$-th root, if $2-r<0$). Then we can compute again
\[
 \det(g_1)\det(g_2) = \lambda^r\det(g_2)=\left(\det(g_2)^{\frac{1}{2-r}}\right)^r\cdot\left(\det(g_2)^{\frac{1}{2-r}}\right)^{2-r},
\]
where we choose the same root of $\det(g_2)\in\C^*$ in the second factor as we did in the definition of $\lambda$. Then by the ordinary power law $z^mz^n=z^{m+n}$ for $z\in\C$ and $m,n\in\Z$ we get
\[
 \det(g_1)\det(g_2)=(\det(g_2)^{\frac{1}{2-r}})^2=\lambda^2
\]
and from this as before that $g_2$ is in the image of $\pi$.
\bigskip

\noindent
{\bf 2.2.} Let us compute the isotropy groups in the $6$ representatives of the unimodular Lie algebras and therefore in particular the dimensions of the corresponding orbits $B_i$ with $i\in\{I,II,VI_0,VII_0,VIII,IX\}$.
\bigskip

\begin{description}
 \item[a. Bianchi $I$:] In this case $k=l=0$ and $m=3$ we have obviously $H=G=\GL_3\R$ and $A_0=0$ is the unique fixpoint of the action. Thus the orbit $B_I$ is given by $B_I=\{0\}$ and therefore $\dim B_I=0$.
 \item[b. Bianchi $II$:] Here, as we saw, $g_2\in\GL_2\R$ is arbitrary, $h\in\R^2$ anyway and $g_1\in\R^*$ is determined by $g_2$ through
\[
 g_1^2=g_1\1_1g_1^T=\det(g_1)\det(g_2)\1=g_1\det(g_2),
\]
i.e.,
\[
 g_1=\det(g_2).
\]
We get
\[
 H=\{\left(\begin{array}{cc} \det(g_2) & h \\ 0 & g_2 \end{array}\right):\; g_2\in\GL_2\R,\; h\in\R^2\}
\]
and therefore, as a manifold, $H$ is diffeomorphic to $\GL_2\R\times\R^2$ which is $6$-dimensional. Thus the orbit $B_{II}$ is $3$-dimensional.
 \item[c. Bianchi $VI_0$:] Let us denote by $\CO(1,1)$ the {\sl conformal group of the Minkowski plane}, i.e.,
\[
 \CO(1,1):=\{g_1\in\GL_2\R:\; \exists\lambda\in\R^*:\; g_1\1_{1,-1}g_1^T=\lambda\1_{1,-1}\}.
\]
By scaling $g_1$ (with $\sqrt{|\lambda|}$) we can get a $g_0$ satisfying
\[
 g_0\1_{1,-1}g_0^T=\pm\1_{1,-1}
\]
and multiplying elements $g_0$ with $g_0\1_{1,-1}g_0^T=-\1_{1,-1}$ by a fixed reflection $s_0$, e.g., $s_0=\diag(1,-1)$, we get finally an element $\tilde g_0$ satisfying $\tilde g_0\1_{1,-1}g_0^T=\1_{1,-1}$. Therefore the dimension of $\CO(1,1)$ is
\[
 \dim\CO(1,1)=1+\dim\O(1,1)=1+1=2.
\]
Now if $\left(\begin{array}{cc}g_1&h\\0&g_2\end{array}\right)\in H$, then obviously $g_1$ is in $\CO(1,1)$ with $\lambda:=\det(g_1)\det(g_2)$ and it is arbitrary. $\det(g_2)$, and therefore $g_2\in\R^*=\GL_1\R$, is uniquely determined by $g_1$ since $\lambda=\iprod{g_1e_1}{e_1}$ and
\[
 \det(g_2)=\frac{\lambda}{\det(g_1)}\in\{-1,1\}.
\]
Thus we get
\[
 H=\{\left(\begin{array}{cc} g_1 & h \\ 0 & \frac{\iprod{g_1e_1}{e_1}}{\det(g_1)}\end{array}\right)\in\GL_3\R:\; g_1\in\CO(1,1),\; h\in\R^2\}
\]
and therefore $\dim H=2+2=4$. The orbit $B_{VI_0}$ is therefore $5$-dimensional.
 \item[d. Bianchi $VII_0$:] Similarly we introduce the {\sl euclidean conformal group}
\[
 \CO(2)=\{g_1\in\GL_2\R:\; \exists\lambda\in\R^*:\; g_1\1_2g_1^T=\lambda\1_2\}.
\]
In this case, by the way, $\lambda$ can only be positive since $g_1g_1^T$ is positive definite. Similarly as in case c.\ we see that
\[
 H=\{\left(\begin{array}{cc}g_1&h\\0&g_2\end{array}\right):\; g_1\in\CO(2),\; h\in\R^2\},
\]
where we can also set
\[
 g_2=\frac{|\det(g_1)|}{\det(g_1)}=\sgn(\det(g_1))
\]
in this case since $\lambda>0$ and therefore $\lambda=|\det(g_1)|$. Again we have $\dim H=4$ and therefore $\dim B_{VII_0}=5$.
  \item[e. Bianchi $VIII$:] In this case $m=0$ and therefore $\det(g_2)=1$. The condition for $g=g_1$ is now:
\[
 g\1_{2,1}g^T=\det(g)\1_{2,1}.
\]
It follows (as in (\ref{eqn:Det}))
\[
 \det(g)^2(-1)=\det(g)^3(-1),
\]
i.e.,
\[
 \det(g)=1.
\]
This shows
\begin{eqnarray*}
 H & = & \{g\in\GL_3\R:\; g\1_{2,1}g^T=\1_{2,1}\; \mbox{und}\;\det(g)=1\} \\
   & = & \SO(2,1)
\end{eqnarray*}
which is $3$-dimensional. The corresponding orbit $B_{VIII}$ is therefore $6$-dimensional.
  \item[f. Bianchi $IX$:] Similarly to the case e.\ one sees that $H=\SO(3)$ and therefore the orbit $B_{IX}$ is also $6$-dimensional.
 \end{description}

\noindent
Our computations imply now a complete description of the orbit structure in the $6$-dimensional $\GL_3$-representation $W^\prime=\Lambda^2(\R^3)^*\otimes S^2\R^3$: All orbits $B$ are (double-) cones, i.e., if $C\in B$ and $\lambda\in\R^*$, then $\lambda C\in B$ as well. There exist two (big) open orbits $B_{VIII}$ and $B_{IX}$ which are separated by a $5$-dimensional cone. This $5$-dimensional cone has itself two (big) open orbits $B_{VI_0}$ and $B_{VII_0}$ which are separated by a $3$-dimensional cone $L$ and this one consists finally of only two orbits, namely its cusp $B_I$ and the rest $B_{II}$.

An orbit is in the closure of another orbit if and only if it has lower dimension.

It is now also easy to see that $W^\prime$ is irreducible. If an invariant subspace $U\sub W^\prime$ contains the orbit $B_{VIII}$ or $B_{IX}$, then $U=W^\prime$ since $U$ contains an open non empty set, of course. If $U$ contains, let's say, the orbit $B_{II}$, then it contains the matrices $\diag(1,0,0)$, $\diag(0,1,0)$ and $\diag(0,0,1)$ and therefore all diagonal matrices. But every orbit has elements which are diagonal, from which we see again that $U$ must be the whole space. Similarly one can argue with the orbits $B_{VI_0}$ and $B_{VII_0}$.

The complex case is handled quite similarly, in fact it is easier. The orbit of the simple Lie algebra has dimension $6$ and is open and dense. The orbit of the solvable Liealgebra has dimension $5$, the orbit of the nilpotent one has dimension $3$ and the abelian has dimension $0$. The closure of an orbit contains all orbits of lower dimension and the representation is irreducible.
\bigskip

\noindent
{\bf 2.3.} Let us classify now also the isomorphism classes for the remaining $B$-classes, i.e., the non unimodular $3$-dimensional Lie algebras. For this we want to prove first a simple lemma which we want to formulate just in the context of group actions on sets.

\begin{lemma}
Suppose $M$ and $N$ are sets on which a group $G$ is acting by permutations. Let further be $f\colon M\to N$ be a $G$-equivariant map and $y_0\in N$. Finally we assume that the action of $G$ on $N$ is transitive. Let us now denote by $F\sub M$ the fiber of $f$ over $y_0$, $F=f^{-1}(y_0)$, and let $H\sub G$ be the stabilizer group in $G$ of $y_0$. Then $F$ is, of course, $H$-invariant and we denote by $\pi_F\colon F\to F/H=:\bar F$ resp.\ $\pi_M\colon M\to M/G=:\bar M$ the natural projections onto the orbit spaces of $H$ resp.\ $G$. Then the inclusion $i\colon F\to M$ induces a unique map $\bar i\colon\bar F\to \bar M$ with $\pi_M\circ i=\bar i\circ\pi_F$. This map $\bar i$ is bijective.
\end{lemma}

\noindent
\begin{proof}
 \begin{enumerate}
  \item The fiber $F$ is in fact $H$-invariant since for $x\in F$ and $h \in H$ and the equivariance of $f$ we have
\[
 f(h.x)=h.f(x)=h.y_0=y_0.
\]
If $x_1,x_2\in F$ are in the same $H$-orbit, i.e., $x_2=h.x_1$ for some $h\in H$, then $i(x_1),i(x_2)\in M$ are in the same $G$-orbit, of course, since $i$ is $H$-equivariant, if you want,
\[
 i(x_2)=i(h.x_1)=h.i(x_1).
\]
Therefore $\pi_M\circ i(x_1)=\pi_M\circ i(x_2)$ and from the universal property of the quotient $\pi_F\colon F\to\bar F$ we get a unique map $\bar i\colon\bar F\to\bar M$ with $\bar i\circ\pi_F=\pi_M\circ i$.

  \item $\bar i$ is surjective: let $c=[x]_M\in \bar M$ be arbitrary. Then there exists a representative $x\in F$ for $c$, because $G$ acts transitively on $N$. Indeed, if we let first $x^\prime\in M$ be an arbitrary representative of $c$ in $M$, we choose $g\in G$ with $g.f(x^\prime)=y_0$. Then $x:=g.x^\prime\in F$ and another representative of $c$. Finally we see that
\[
 \bar i([x]_F)=\bar i\circ\pi_F(x)=\pi_M\circ i(x)=[x]_M=c.
\]

   \item $\bar i$ is also injective: So let $\bar i([x_1])=\bar i([x_2])$. Then 
\[
 [x_1]_M=\pi_M\circ i(x_1)=\bar i\circ\pi_F(x_1)=\bar i([x_1])=\bar i([x_2])=\ldots=[x_2]_M,
\]
and we can find a $g\in G$ so that $g.x_1=x_2$. But
\[
 g.y_0=g.f(x_1)=f(g.x_1)=f(x_2)=y_0,
\]
i.e., $g\in H$ and therefore $[x_1]_F=[x_2]_F$.
 \end{enumerate}
\end{proof}

\noindent
We consider now the natural $G$-equivariant projection $p\colon X\to W^\prime$ of our Lie variety $X\sub\Lambda^2V^*\otimes V$ onto $W^\prime=\Lambda^3V^*\otimes S^2V$ (with $G=\Aut(V)$ as always). Then we want to apply the lemma to the restriction of $p$ to the preimages $p^{-1}(B_i)$, where $B_i\sub W^\prime$ is one of the $6$ (resp.\ $4$) unimodular orbits of $G$ on $W^\prime$,
\[
 p_i\colon p^{-1}(B_i)\to B_i
\]
($i=I,II,VI_0,VII_0,VIII,IX$). As our reference points $y_i\in B_i$ we choose our representatives $M_i$ coming from our choice of matrices $A_i\in\Sym_2(\R^3)$. Since we already computed the stabilizer groups $H_i\sub G$ and also the fibers $F_i\sub X$ it remains to study the orbits of $H_i$ on $F_i$. But $H_i$ acts on $F_i=\{M_i\}\times\Deg(M_i)\sub W^\prime\oplus V^*\cong\Lambda^2V^*\otimes V$ by
\[
 h.(M_i,\alpha)=(M_i,h.\alpha),
\]
since $h.M_i=M_i$. Here the action of $H_i\sub G$ on $\Deg(M_i)\sub V^*$ is the restriction to $H_i$ on $\Deg(M_i)$ of the natural action of $G=\Aut(V)$ on $V^*$ (cf.\ 2.1). In coordinates this comes down to the natural action of $\pi(H)\sub\GL_m\R$ on $\Deg(M)\cong\R^m$. For the case where $\pi(H)=\GL_m\R$, there exist therefore only two orbits, namely the origin and the rest. Thus we get the following cases:

 \begin{description}
  \item[Bianchi $I$:] In this case $M=0$, i.e., $m=3$ and $\Deg(M)=V^*$ and over $M_I$ there is only one other orbit, namely $j(V^*\setminus\{0\})\sub W_2$ as we essentially knew before. This type is called Bianchi $V$.

  \item[Bianchi $II$:] In this case we have $k=1$, $l=0$ and $m=2$ and therefore again, there is only one orbit over $M_{II}$ which is usually called Biachi $IV$.

  \item[Bianchi $VIII$ and $IX$:] In this case $M$ is non degenerate, i.e., $m=0$ and there is no orbit at all over $M$ in this case.

  \item[Bianchi $VI_0$ and $VII_0$:] In this case $(k,l)=(1,1)$ resp.\ $(k,l)=(2,0)$ and therefore $m=1$. From corollary \ref{cor:} we get $\pi(H)=\{\pm 1\}\sub\R^*=\GL_1\R$, i.e., for every $\pm h$ we get a new orbit. For the Bianchi $VI$ case one usually chooses $h\leq 0$ and for Bianchi $VII$ one chooses $h\geq 0$.
 \end{description}

\noindent
This classifies all real Lie algebras of dimension $3$.

The complex case works similarly. Over $\Bia(VIII)$ and $\Bia(IX)$ there is nothing and over $\Bia(I)$ and $\Bia(II)$ there is exactly one other orbit (which we also call $\Bia(V)$ resp.\ $\Bia(IV)$), and over $\Bia(VII_0)$ there is a series $\Bia(VII_h)$ for $h\in\C^*$ and $\Bia(VII_h)$ is the same Lie algebra as $\Bia(VII_{-h})$ (and only these are isomorphic).

We want to give a choice of structure constants $(C_{ij}^k)$ (for $1\leq i<j\leq 3$ and $1\leq k\leq 3$)  for each structure given by a choice of a representative $C_i\in B_i$ (and a choice of vector space $V$ and a basis $(e_i)$ of $V$, where we choose $V=\R^3$ and $(e_i)$ the canonical basis of $\R^3$, of course). For this observe first that the Lie structure $C\in\Lambda^2V^*\otimes V$ given in the $A$-class by an element
\[
 M\in\Lambda^3V^*\otimes S^2V\cong\Sym_2(V^*,\Lambda^3V^*)
\]
($V=\R^3$) by the diagonal matrix $A=\diag(\mu_1,\mu_2,\mu_3)$ is given on the canonical basis $(e_i)$ by
\begin{equation} \label{eqn:Lie}
 \Lie{e_i}{e_j}=\sum_{k=1}^3\varepsilon_{ijk}\mu_ke_k,
\end{equation}
where $\varepsilon_{ijk}=1$, if $(i,j,k)$ is an even permutation of $(1,2,3)$, $\varepsilon_{ijk}=-1$, if $(i,j,k)$ is an odd permutation and $\varepsilon_{ijk}=0$ otherwise. In fact, for the dual basis $(\lambda^i)$ we have that
\begin{eqnarray*}
 \iprod{\lambda^l}{\Lie{e_i}{e_j}} & = & \iprod{\lambda^l}{C(e_i,e_j)}=\iprod{\lambda^l}{\Tr(M)(e_i,e_j)} \\ & = & M(e_i,e_j,e_k)(\lambda^k,\lambda^l)=\varepsilon_{ijk}M_{123}(\lambda^k,\lambda^l)=\varepsilon_{ijl}\mu^l,
\end{eqnarray*}
i.e.,
\[
 \Lie{e_i}{e_j}=\sum_k\varepsilon_{ijk}\mu_ke_k.
\]

On the other hand the structure constants $C_{ij}^k=\iprod{\lambda^k}{\Lie{e_i}{e_j}}$ for the structures $j(\nu)$ with $\nu\in V^*$ are given by
\begin{eqnarray*}
 C_{ij}^k & = & \iprod{\lambda^k}{j(\nu)(e_i,e_j)}=\iprod{\lambda^k}{\nu(e_i)e_j-\nu(e_j)e_i} \\ & = & \nu(e_i)\delta_j^k-\nu(e_j)\delta_i^k.
\end{eqnarray*}
In the special case $\nu=h\lambda^3$ this comes down to
\begin{equation} \label{eqn:Lie2}
 \Lie{e_1}{e_2}=0,\quad\Lie{e_2}{e_3}=-he_2,\quad\Lie{e_3}{e_1}=he_1.
\end{equation}
For the $4$ new $B$-classes we choose as representatives $(M,\nu)$ with $\nu=\lambda^3$ for the Bianchi $IV$ and $V$ case and $\nu=h\lambda^3$ (with $h<0$ for $\Bia(VI_h)$ and $h>0$ for $\Bia(VII_h$)) and $M$ the choice from above for the corresponding $A$-class. In this way we get:

\begin{theorem} (Bianchi). Let $L$ be a real $3$-dimensional Lie algebra. Then $L$ is isomorphic to exactly one of the following Lie algebras $(\R^3,\Lie{\cdot}{\cdot})$, where the Lie bracket is given on the canonical basis as follows:
\begin{description}
 \item[Bianchi $I$:] $\Lie{e_1}{e_2}=\Lie{e_2}{e_3}=\Lie{e_3}{e_1}=0$
 \item[Bianchi $II$:] $\Lie{e_1}{e_2}=0,\;\Lie{e_2}{e_3}=e_1,\;\Lie{e_3}{e_1}=0$
 \item[Bianchi $IV$:] $\Lie{e_1}{e_2}=0,\;\Lie{e_2}{e_3}=e_1-e_2,\;\Lie{e_3}{e_1}=e_1$
 \item[Bianchi $V$:] $\Lie{e_1}{e_2}=0,\;\Lie{e_2}{e_3}=e_2,\;\Lie{e_3}{e_1}=e_1$
 \item[Bianchi $VI_h$ ($h\leq 0$):] $\Lie{e_1}{e_2}=0,\;,\Lie{e_2}{e_3}=e_1-he_2,\;\Lie{e_3}{e_1}=he_1-e_2$
 \item[Bianchi $VII_h$ ($h\geq 0$):] $\Lie{e_1}{e_2}=0,\;\Lie{e_2}{e_3}=e_1-he_2,\;\Lie{e_3}{e_1}=he_1+e_2$
 \item[Bianchi $VIII$:] $\Lie{e_1}{e_2}=-e_3,\;\Lie{e_2}{e_3}=e_1,\;\Lie{e_3}{e_1}=e_2$
 \item[Bianchi $IX$:] $\Lie{e_1}{e_2}=e_3,\;\Lie{e_2}{e_3}=e_1,\;\Lie{e_3}{e_1}=e_2$
\end{description}
\end{theorem}

\noindent
\begin{proof} The structure constants in the theorem are given directly from the summation of the structure constants of the Bianchi $A$ classes (\ref{eqn:Lie}) and the $B$-class $j(\nu)$ in (\ref{eqn:Lie2}).
\end{proof}
\bigskip

\noindent
The Bianchi $III$ class which the reader might miss at this point, or even earlier, coincides with $\Bia(VI_{-1})$. $L=\Bia(VI_{-1})$ is in fact somehow special inside the series $(\Bia(VI_h))_{h<0}$, since it is the only of these Lie algebras which has a only $1$-dimensional {\sl derivative} $L^\prime:=\Lie{L}{L}$, the others have a $2$-dimensional derivative. Even among all Lie algebras there is only one other Lie algebra with a $1$-dimensional derivative, namely the {\sl Heisenberg Lie algebra} $\Bia(II)$. The natural exact sequence 
\[
 0\to L^\prime\to L\to L/L^\prime\to 0
\]
gives another approach to classify all solvable $3$-dimensional Lie algebras (cf. \cite{KKL}).

$L=\Bia(III)$ is, by the way, also the only non-trivial product of lower dimensional Lie algebras, namely $L\cong\R\times h$, where $h$ denotes the non abelian $2$-dimensional Lie algebra.
\bigskip

\noindent
{\bf 2.4} Let us finally complete our list of stabilizer groups and therefore the automorphism groups of the $B$-classes $C\in X$. Of course, these are subgroups of the isotropy groups of the corresponding $A$-class $\frac{1}{2}\Tr\circ p(C)\in W_1$.
 \begin{description}
  \item[a. Bianchi $V$:] As above we choose the representative $\nu=\lambda^3$ and $p(C)=0$, of course. From that it follows that the isotropy group $H=G_V$ is given by
\[
 H=\{g\in\GL_3\R:\; g=\left(\begin{array}{cc} g_1 & 0 \\ h & 1 \end{array}\right)\; \mbox{with $g_1\in\GL_2\R$ and $h\in\R^2$}\},
\]
which is isomorphic to the {\sl affine group} $\Aff_2(\R)\cong\GL_2\R)\times_\sigma\R^2$ (semidirect product w.r.t. the obvious representation $\sigma$ of $\GL_2\R$ on $\R^2$). The orbit dimension is therefore $9-6=3$, which we have known before. The closure of that orbit contains (besides itself) only the orbit $B_I$ (the origin). Of course, the closure is just the irreducible submodule $W_2\cong V^*$.

  \item[b. Bianchi $IV$:] The stabilizer group $H$ for $(M,\nu)=(\diag(1,0,0),\lambda^3)$ is given by
\[
 H=\{\left(\begin{array}{ccc} \mu & h_1 & h_2 \\ 0 & \mu & 0 \\ 0 & h_3 & 1 \end{array}\right)\in\GL_3\R:\; \mu\in\R^*,\; h_1,h_2,h_3\in\R\}.
\]
This has dimension $4$ and therefore the corresponding orbit $B_{IV}$ has dimension $5$, which we also knew, since $B_{II}$ has dimension $3$ and $p\colon B_{IV}\to B_{II}$ has $2$-dimensional fibers which are affine subspaces of dimension $2$. In the closure of $B_{IV}$ there are $B_{IV}$ (of course), $B_{II}$, $B_V$ and $B_I$ which is quite clear.

  \item[c. Bianchi $VI_h$ ($h<0$):] The stabilizer groups in $(M,\nu)=(\diag(1,-1,0),$ $h\lambda^3)$  are for all $h<0$ the same and are of index $2$ in the stabilzer group of $\Bia(VI_0)$. It is given by
\[
 H=\{\left(\begin{array}{cc} g_1 & h \\ 0 & 1 \end{array}\right)\in\GL_3\R:\; g_1\in\CO(1,1),\; h\in\R^2\},
\]
i.e., this is again a semidirect product $\CO(1,1)\times_\sigma\R^2$. The orbits have all dimension $5$ which we knew before, since the projection $p\colon B_{VI_h}\to B_{VI_0}$ is a $2:1$-cover over the $5$-dimensional $B_{VI_0}$. In the closure of $B_{VI_h}$ we find the orbits $B_{VI_h}$, $B_{II}$ and $B_I$ (and not $B_{IV}$ or $B_V$). However, $B_{IV}$ is in the closure of the union $\bigcup_{h<0}B_{VI_h}$ as is $B_V$ and $B_{VI_0}$.

  \item[d. Bianchi $VII_h$ ($h>0$):] This case is essentially the same as the former. $B_{VII_h}$ is $5$-dimensional, it lies $2:1$ over $B_{VII_0}$ and has in its closure (besides itself) the orbits $B_{II}$ and $B_I$. The union $\bigcup_{h>0}B_{VII_h}$ is $6$-dimensional and has additionally the orbits $B_{IV}$, $B_V$ and $B_{VII_0}$ in its closure. The stabilizer group is
\[
 H=\{\left(\begin{array}{cc} g_1 & h \\ 0 &1 \end{array}\right)\in\GL_3\R:\; g_1\in\CO(2),\; h\in\R^2\}\cong\CO(2)\times_\sigma\R^2.
\]
 \end{description}

\noindent
We let it to the reader to determine the stabilizer groups of the $B$-classes in the complex case as well as the orbit dimensions and the topology of the orbit space $X/G$.

\end{document}